\documentclass[12pt]{article}
\usepackage{graphicx}
\usepackage{amsmath}

\newcommand{\bb}{\begin{eqnarray*}}

\newcommand{\beq}{\begin{equation}}
\newcommand{\eeq}{\end{equation}}

\newcommand{\ee}{\end{eqnarray*}}

\newtheorem{theorem}{Theorem}

\newtheorem{corollary}[theorem]{Corollary}

\newtheorem{definition}[theorem]{Definition}

\newtheorem{lemma}[theorem]{Lemma}

\newtheorem{proposition}[theorem]{Proposition}
\newtheorem{remark}[theorem]{Remark}
\newtheorem{comment}[theorem]{Comment}

\input{tcilatex}
\begin{document}

\begin{center}
\textbf{\Large Functional moderate deviations for triangular arrays and
applications}\vskip15pt

Florence Merlev\`{e}de $^{a}$ \textit{and\/} Magda Peligrad $^{b}$ \footnote{%
Supported in part by a Charles Phelps Taft Memorial Fund grant and NSA\
grant, H98230-07-1-0016.}
\end{center}

$^{a}$ Universit\'{e} Paris 6, LPMA and C.N.R.S UMR 7599, 175 rue du
Chevaleret, 75013 Paris, FRANCE

\vskip10pt $^{b}$ Department of Mathematical Sciences, University of
Cincinnati, PO Box 210025, Cincinnati, Oh 45221-0025 \vskip10pt



\begin{center}
\textbf{Abstract}\vskip10pt
\end{center}

\quad Motivated by the study of dependent random variables by coupling with
independent blocks of variables, we obtain first sufficient conditions for
the moderate deviation principle in its functional form for triangular
arrays of independent random variables. Under some regularity assumptions
our conditions are also necessary in the stationary case. The results are
then applied to derive moderate deviation principles for linear processes,
kernel estimators of a density and some classes of dependent random
variables.

\section{Introduction}

\qquad In recent years substantial progress was achieved in obtaining
necessary and sufficient condition for the moderate deviations behavior of
sums of independent identically distributed random variables. Papers by
Ledoux (1992) and Arcones (2003-a, 2003-b, 2003-c), among others, are steps
in this direction. These works show that the moderate deviation principle
can be applied for i.i.d. sequences, even when the moment generating
function is not defined in a neighborhood of zero. Due to its invariant
nature, a natural question is to treat triangular arrays of random
variables. Some sufficient conditions for bounded triangular arrays are
contained in Lemma 2.3 in Arcones (2003-a) and also in the results by
Puhalskii (1994) about triangular arrays of martingale differences. Djellout
(2002) studied this problem for not necessarily stationary martingale
differences sequences.

\qquad In this paper we derive sufficient conditions for the moderate
deviation principle in its functional form for triangular arrays of
independent random variables. In the stationary case and under some
regularity conditions, the condition is necessary as well. These results
open the way to address the moderate deviation principle for classes of
dependent random variables that were not studied so far, by dividing the
variables in blocks that are further approximated by a triangular array of
independent random variables. As a matter of fact this was the initial
motivation of our study. The results are used to treat general linear
processes, Kernel estimators of a density, and some dependent structures
including classes of strong mixing sequences.

\quad The moderate deviation principle is an intermediate estimation between
central limit theorem and large deviation. We shall assume for the moment
that we have a triangular array of independent, centered and square
integrable random variables $(X_{n1},$ $X_{n2},...,X_{nk_{n}})$, where $%
k_{n} $ is a sequence of integers. Denote by 
\begin{equation*}
\text{$S_{n,0}=0$, $S_{nl}=\sum_{j=1}^{l}X_{nj}$, $S_{n}=%
\sum_{j=1}^{k_{n}}X_{nj}$, $\sigma _{nj}^{2}=\mathrm{Var}(X_{nj})$, $%
s_{n}^{2}=\sum_{j=1}^{k_{n}}\sigma _{nj}^{2}$ and $s_{ni}^{2}=\sum_{j=1}^{i}%
\sigma _{nj}^{2}$}.
\end{equation*}%
In the rest of the paper MDP stays for Moderate Deviation Principle.

\begin{definition}
We say that the $MDP$ holds for $s_n^{-1}S_{n}$ with the speed $a_{n}
\rightarrow 0$ and rate function $I(t)$ if for each $A$ Borelian, 
\begin{eqnarray}
-\inf_{t\in A^{o}}I(t) & \leq& \lim\inf_{n}a_{n}\log {\mathbf{P}}(\frac{%
\sqrt{a_{n}}}{s_n}S_{n} \in A)  \notag \\
& \leq & \lim\sup_{n}a_{n}\log {\mathbf{P}}(\frac{\sqrt{a_{n}}}{s_n}S_{n}\in
A)\leq-\inf_{t\in\bar{A}}I(t) \, .  \label{MDPCLT}
\end{eqnarray}
\end{definition}

\quad We are also interested to give a more general result concerning the
Donsker process associated to the partial sums.

\begin{definition}
\label{defdonskerline} Let $\{W_{n},n>0\}$ be the family of random variables
on $D[0,1]$ defined as follows:
\end{definition}

$W_{n}(t)=S_{n,i-1}/s_{n}$ for $t\in \lbrack
s_{n,i-1}^{2}/s_{n}^{2},s_{ni}^{2}/s_{n}^{2})$ , $\ $where $1\leq i\leq
k_{n}\ $\ and $W_{n}(1)=S_{n}/s_{n}$.\newline
We say that the family of random variables $\{W_{n},n>0\}$ satisfies the
functional Moderate Deviation Principle (MDP) in $D[0,1]$ endowed with
uniform topology, with speed $a_{n}\rightarrow 0$ and good rate function $%
I(.)$, if the level sets $\{x,I(x)\leq \alpha \}$ are compact for all $%
\alpha <\infty $, and for all Borel sets $\Gamma \in {\mathcal{B}}$ 
\begin{eqnarray}
-\inf_{t\in \Gamma ^{0}}I(t) &\leq &\lim \inf_{n}a_{n}\log {\mathbf{P}}(%
\sqrt{a_{n}}W_{n}\in \Gamma )  \notag \\
&\leq &\lim \sup_{n}a_{n}\log {\mathbf{P}}(\sqrt{a_{n}}W_{n}\in \Gamma )\leq
-\inf_{t\in \bar{\Gamma}}I(t)\,.  \label{mdpkey}
\end{eqnarray}

\quad Our first result is:

\begin{theorem}
\label{keyresult}Assume that $(X_{n1},$ $X_{n2},...,X_{nk_{n}})$ is a
triangular array of independent centered and square integrable random
variables. Assume $a_{n}\rightarrow 0$ and that for any $\beta >0$%
\begin{equation}
\lim \sup_{n\rightarrow \infty }a_{n}\sum_{j=1}^{k_{n}}{\mathbf{E}}([\exp
\beta \frac{|X_{nj}|}{\sqrt{a_{n}}s_{n}}]I(\sqrt{a_{n}}s_{n}<|X_{nj}|<s_{n}/%
\sqrt{a_{n}})=0\,,  \label{lindebergmain2}
\end{equation}%
\begin{equation}
\lim \sup_{n}a_{n}\log {\mathbf{P}}(\max_{1\leq j\leq k_{n}}|X_{nj}|\geq
s_{n}/\sqrt{a_{n}})=-\infty  \label{neg2}
\end{equation}%
and for any $\epsilon >0$ 
\begin{equation}
\frac{1}{s_{n}^{2}}\sum_{j=1}^{k_{n}}{\mathbf{E}}[X_{nj}^{^{2}}I(|X_{nj}|%
\geq \epsilon s_{n}\sqrt{a_{n}})]\rightarrow 0\text{ }.  \label{lindeberg2}
\end{equation}%
Then $\{W_{n},n>0\}$ satisfies MDP in $D[0,1]$ with speed $a_{n}$ and rate
function $I(.)$ defined by 
\begin{align}
I(z)& =\frac{1}{2}\int_{0}^{1}(z^{\prime }(u))^{2}du\,\text{if }z(0)=0\text{
\ and }z\text{ is absolutely continuous}  \label{rate} \\
& \text{ and }\infty \text{ otherwise.}  \notag
\end{align}
\end{theorem}

\begin{comment}
\label{commentcond copy(1)}Under the assumptions of the theorem, we have in
particular that $\{s_{n}^{-1}\sum_{j=1}^{k_{n}}X_{nj}\}$ satisfies the MDP
with speed $a_{n}$ and rate $I(t)=t^{2}/2$.
\end{comment}

\quad Standard computations show that all the conditions of Theorem \ref%
{keyresult} are satisfied if we impose the unique condition (\ref{lindmain21}%
) below, that can be viewed as a generalized Lindeberg's condition. So we
can state:

\begin{corollary}
\label{one(1)}Assume $(X_{n1},$ $X_{n2},...,X_{nk_{n}})$ is a triangular
array of independent centered and square integrable random variables. Assume 
$a_{n}\rightarrow 0$, and for any $\epsilon >0$ and any $\beta >0$, 
\begin{equation}
\lim \sup_{n\rightarrow \infty }a_{n}\sum_{j=1}^{k_{n}}{\mathbf{E}}([\exp
\beta \frac{|X_{nj}|}{\sqrt{a_{n}}s_{n}}]I(|X_{nj}|>\epsilon \sqrt{a_{n}}%
s_{n})=0 \, .  \label{lindmain21}
\end{equation}%
Then the conclusion of Theorem \ref{keyresult} is satisfied.
\end{corollary}

\quad Simple computations involving Chebyshev's inequality and integration
by parts (see Appendix) lead to conditions imposed to the tails
distributions of the random variables involved.

\begin{comment}
\label{comment}Condition (\ref{lindebergmain2}) is equivalent to : There is
a constant $C_{1}$ with the following property: for any $\beta >0$ there is $%
N(\beta )$ such that for $n>N(\beta )$ 
\begin{equation}
a_{n}\sum_{j=1}^{k_{n}}{\mathbf{P}}(|X_{nj}|>u\sqrt{a_{n}}s_{n})\leq
C_{1}\exp (-\beta u)\text{ for all }1\leq u\leq 1/a_{n}\,.  \label{vv}
\end{equation}%
Condition (\ref{lindmain21}) is equivalent to : There is a constant $C_{1}$
with the property that for any $\epsilon >0$ and any $\beta >0$, there is $%
N(\epsilon ,\beta )$ such that for $n>N(\epsilon ,\beta )$, the inequality
in relation (\ref{vv}) is satisfied for all $u\geq \epsilon $.
\end{comment}

\medskip \quad If we impose some regularity assumptions the conditions
simplify.

\medskip

\noindent \textbf{RC} The functions $f(n)=s_{n}^{2}a_{n}$ and $%
g(n)=s_{n}^{2}/a_{n}$ are strictly increasing to infinite, and the function $%
l(n)=s_{n}^{2}/k_{n}$ is nondecreasing.

\medskip \quad Assuming \textbf{RC}, we construct the strictly increasing
continuous function $f(x)$ that is formed by the line segments from $%
(n,f(n)) $ to $(n+1,f(n+1))$. Similarly we define $g(x)$ and denote by $%
c(x)=f^{-1}(g(x)).$

\begin{corollary}
\label{onenonbis} Assume $(X_{n1},$ $X_{n2},...,X_{nk_{n}})$ is a triangular
array of independent, centered and square integrable random variables.
Assume $a_{n}\rightarrow 0$, the regularity conditions \textrm{\textbf{RC}}
hold and 
\begin{equation}
a_{n}\log (\sup_{n\leq m\leq c(n+1)}\sup_{1\leq i\leq k_{m}}k_{n}{\mathbf{P}}%
(|X_{mi}|>s_{n}/\sqrt{a_{n}}))\rightarrow -\infty \text{ as }n\rightarrow
\infty \,.  \label{onecondm}
\end{equation}%
Assume in addition that (\ref{lindeberg2}) is satisfied. Then the conclusion
of Theorem \ref{keyresult} holds.
\end{corollary}

In the sequel we shall denote by $[x]$ the integer part of $x.$

\begin{remark}
\label{station} In the case where $(X_{n})_{n\geq 0}$ is a sequence of
i.i.d. r.v's with mean zero and finite second moment the conditions of
corollary \ref{onenonbis} simplify. If $a_{n}\searrow 0$, $na_{n}\nearrow
\infty $ and%
\begin{equation}
a_{n}\log n{\mathbf{P}}\Big (|X_{0}|>\frac{\sigma \sqrt{n}}{\sqrt{a_{n}}}%
\Big )\rightarrow -\infty  \label{equiarcones}
\end{equation}%
then, the conclusion of Theorem 3 holds with $W_{n}(t)=n^{-1/2}%
\sum_{j=1}^{[nt]}X_{j}$. Moreover, condition (\ref{equiarcones}) is
necessary for the moderate deviation principle in this case. This result for
i.i.d. is contained in Arcones (Theorem 2.4, 2003-a).
\end{remark}

\quad For the sake of applications we give a sufficient condition in terms
of the moments of $X_{n,i}$.

\begin{proposition}
\label{forappl} Assume $(X_{n1},$ $X_{n2},...,X_{nk_{n}})\ $is a triangular
array of independent centered and square integrable random variables. Assume
that there exists $n_{0}$ such that for each $n\geq n_{0}$ and $1\leq k\leq
k_{n}$ there are nonnegative numbers $A_{nk}$ and $B_{n}$ such that for each 
$m\geq 3$ 
\begin{equation}
{\mathbf{E}}|X_{n,k}|^{m}\leq m!A_{nk}^{m}B_{n}\,.  \label{mom}
\end{equation}%
Assume in addition that $a_{n}\rightarrow 0$, 
\begin{equation}
A_{n,k}=o(\sqrt{a_{n}}s_{n})\text{ as }n\rightarrow \infty \text{ uniformly
in }k  \label{condAk1}
\end{equation}%
and there is a positive constant $C$ such that 
\begin{equation}
\frac{B_{n}}{s_{n}^{2}}\sum_{j=1}^{k_{n}}|A_{nj}|^{2}\leq C\text{ for all }%
n\geq n_{0}\,.  \label{condAk2}
\end{equation}%
Then the conclusion of Theorem \ref{keyresult} holds.
\end{proposition}

\section{Applications}

\subsection{A class of Linear processes}

\qquad In this section, we consider a sequence $\{\xi _{k}\}_{k\in \mathbf{Z}%
}$ of i.i.d. and centered random variables such that ${\mathbf{E}} (\xi
_{0})^{2}=\sigma ^{2}>0$ and let $\{c_{ni},1\leq i\leq k_{n}\}$ be a
triangular array of numbers. Many statistical procedures produce estimators
of the type 
\begin{equation}
S_{n}=\sum_{i=1}^{k_{n}}c_{ni}\xi _{i}\,.  \label{deflingen}
\end{equation}%
\qquad For instance, consider the fixed design regression problem $%
Z_{k}=\theta q_{k}+\xi _{k}$, where the fixed design points are of the form $%
q_{k}=1/g(k/n)$ where $g(.)$ is a function. To analyze the error of the
estimator $\hat{\theta}=n^{-1}\sum_{k=1}^{n}Z_{k}g(k/n)$, we are led to
study the behavior of processes of the form (\ref{deflingen}). \newline
\quad Setting 
\begin{equation}
s_{n}^{2}=\sigma ^{2}\sum_{i=1}^{k_{n}}c_{ni}^{2}\,,  \label{defsigmalingen}
\end{equation}%
we are interested to give sufficient conditions for the moderate deviation
principle for $\displaystyle\frac{S_{n}}{s_{n}}$ and also for the stochastic
process $W_{n}(\cdot )$ defined in Definition {\ref{defdonskerline}} with $%
X_{n,i}=c_{ni}\xi _{i},$ $s_{ni}^{2}=\sigma ^{2}\sum_{j=1}^{i}c_{nj}^{2}$.
In order for the Lindeberg's condition (\ref{lindeberg2}) to be satisfied we
shall impose the following condition%
\begin{equation}
\frac{1}{\sqrt{a_{n}}s_{n}}\max_{1\leq j\leq k_{n}}|c_{nj}|\rightarrow 0\,\
\ {\text{as }}\,n\rightarrow \infty \,.  \label{C1lingen1}
\end{equation}

\quad By applying Corollary \ref{onenonbis} we easily obtain the following
result

\begin{proposition}
\label{proplingen} Let $S_{n}$ and $s_{n}^{2}$ be defined by (\ref{deflingen}%
) and (\ref{defsigmalingen}). Assume that $a_{n}\rightarrow 0$, condition (%
\ref{C1lingen1}) holds and the regularity conditions \textrm{\textbf{RC}}.
Denote by $C_{n}=\sup_{n\leq m\leq c(n+1)}\sup_{1\leq i\leq k_{m}}|c_{m,i}|$
and assume that the following condition holds 
\begin{equation}
a_{n}\log (k_{n}{\mathbf{P}}(|\xi _{0}|>s_{n}/C_{n}\sqrt{a_{n}}))\rightarrow
-\infty \text{ as }n\rightarrow \infty \,.  \label{C1lingen}
\end{equation}%
Then $\{W_{n}(\cdot )\}$ satisfies the MDP in $D[0,1]$ with speed $a_{n}$
and rate $I(\cdot )$ defined in Theorem \ref{keyresult}.
\end{proposition}

\quad Notice that the variable $\xi _{0}$ is not required to have moment
generating functions. As a matter of fact, by using Proposition \ref{forappl}
we can easily derive.

\begin{proposition}
\label{proplingen1} Let $S_{n}$ and $s_{n}^{2}$ be defined by (\ref%
{deflingen}) and (\ref{defsigmalingen}). Assume that $a_{n}\rightarrow 0$
and condition (\ref{C1lingen1}) holds. Assume that for some positive
constant $K,$ 
\begin{equation}
{\mathbf{E}}(|\xi _{0}|^{m})\leq m!K^{m}\,\text{\ for all }m\in {\mathbf{N}}%
\,.  \label{momcond}
\end{equation}%
Then $\{W_{n}(\cdot )\}$ satisfies the MDP in $D[0,1]$ with speed $a_{n}$
and rate $I(\cdot )$ defined in Theorem \ref{keyresult}.
\end{proposition}

\quad To give a few examples, notice that if the double sequence $%
|c_{m,i}|_{m,i}$ is uniformly bounded by a constant, condition (\ref%
{C1lingen1}) is verified provided $\lim_{n\rightarrow \infty
}a_{n}s_{n}^{2}=\infty .$ Moreover, if for each $n$ fixed, the sequence $%
\{|c_{nj}|\}_{j\geq 1}$ is increasing and satisfies the regularity
assumption $\sum_{i=1}^{n}c_{ni}^{2}\sim nc_{n,n}^{2}$, then condition (\ref%
{C1lingen1}) is satisfied if $na_{n}\rightarrow \infty $. This is the case
for instance when $c_{ni}^{2}=c_{i}^{2}=h(i)$ with $h(x)$ a slowly varying
increasing function.

\quad Of course, if smaller classes of random variables $(\xi _{k})_{k\in {%
\mathbf{Z}}}$ are considered, such as bounded or sub-gaussian variables, a
requirement weaker then (\ref{C1lingen1}) may guaranty MDP. We give here an
example showing that condition (\ref{C1lingen1}) of Proposition \ref%
{proplingen1} is necessary when the random variables $(\xi _{k})_{k\in {%
\mathbf{Z}}}$ satisfy only a condition of type (\ref{momcond}).

\quad Assume $(X_{i},i\in {\mathbf{Z}})$ is a sequence of independent
identically distributed random variables with exponential law with mean $1,$
($\mathbf{P}(X_{0}>x)=e^{-x}$), denote $\xi _{n}=X_{n}-1$ and assume the
sequence of constants has the property $\max_{1\leq j\leq n}|c_{nj}|=1$.
Notice first that, for any $t>0$, we have that 
\begin{equation*}
\mathbf{P}\big (\frac{\sqrt{a_{n}}}{s_{n}}|\xi _{0}|\geq t\big )\leq \mathbf{%
P}\big (\max_{1\leq i\leq n}\frac{\sqrt{a_{n}}}{s_{n}}|c_{ni}\xi _{i}|\geq t%
\big )\,.
\end{equation*}%
Also, by standard symmetrization arguments and Levy's inequality (see for
instance Proposition 2.3 in Ledoux and Talagrand (1991)), we get that for
any $t>0$ and $n$ large enough (such that $a_{n}\leq t^{2}/8$), 
\begin{equation*}
{\mathbf{P}}\big (\max_{1\leq i\leq n}\frac{\sqrt{a_{n}}}{s_{n}}|c_{ni}\xi
_{i}|\geq t\big )\leq 2\Big (1-\frac{4a_{n}}{t^{2}}\Big)^{-1}\mathbf{P}\Big (%
\frac{\sqrt{a_{n}}}{s_{n}}|S_{n}|\geq t/2\big )\,.
\end{equation*}%
Now if $\{s_{n}^{-1}S_{n}\}$ satisfies the MDP, then the previous
inequalities entail that necessarily 
\begin{equation*}
\limsup_{n\rightarrow \infty }a_{n}\log \mathbf{P}\big (\frac{\sqrt{a_{n}}}{%
s_{n}}(X_{0}-1)\geq t\big )\leq -\frac{t^{2}}{8}\,.
\end{equation*}%
On an other hand for any $t>0$ 
\begin{equation*}
a_{n}\log \mathbf{P}\big (\frac{\sqrt{a_{n}}}{s_{n}}(X_{0}-1)\geq t\big )%
=-a_{n}-t\sqrt{a_{n}}s_{n}\,.
\end{equation*}%
In order for $\limsup_{n\rightarrow \infty }(-a_{n}-t\sqrt{a_{n}}s_{n})\leq -%
\frac{t^{2}}{8}$ for all $t>0$ we see that necessarily $a_{n}s_{n}^{2}%
\rightarrow \infty $, which implies that condition (\ref{C1lingen1}) is
satisfied since $\max_{1\leq j\leq n}|c_{nj}|=1$.

\medskip

\noindent \textbf{Proof of Proposition \ref{proplingen1}.} We shall apply
Proposition \ref{forappl} with $X_{n,k}=c_{nk}\xi _{k}$. Notice that for all
positive integers $m$, 
\begin{equation*}
\mathbf{E}\big (|X_{nk}|^{m}\big )\leq |c_{nk}|^{m}\mathbf{E}\big (|\xi
_{0}|^{m}\big )\,.
\end{equation*}%
Whence, by using (\ref{momcond}), we get for all positive integers $m$, 
\begin{equation*}
\mathbf{E}\big (|X_{nk}|^{m}\big )\leq m!|c_{nk}|^{m}K^{m}\,.
\end{equation*}%
Then, the conditions of Proposition \ref{forappl} are satisfied and the
result follows.

\subsection{Kernel Estimators of the density}

\qquad In this section we apply our results to obtain a simple MDP in its
functional form for the Kernel estimator at a fix point. Different and
further pointing problems related to the moderate deviation principle for
kernel estimators of the density or of the regression function were
addressed in several papers. For instance, for kernel density estimator,
Louani (1998) addresses the problem of large deviations, Gao (2003) studies
the MDP uniformly in $x$, while Mokkadem, Pelletier and Worms (2005) give
the large and moderate deviation principles for partial derivatives of a
multivariate density. Concerning the kernel estimators of the multivariate
regression, Mokkadem, Pelletier and Thiam (2007) study their large and
moderate deviation principles.

\quad Let $X=(X_{k},k\in {\mathbf{Z)}}$ be a sequence of i.i.d. random
variables. We now impose the following conditions:

\indent {\bf (A.1)} The density function of $X_{0}$ is bounded and
continuous at a fixed point $x$.

\indent {\bf (A.2)} The kernel $K$ is a function such that $\int_{\mathbf{R}%
}K(x)dx=1$ and there exists a positive constant $C$ such that for all
positive integers $m$ 
\begin{equation*}
\int_{\mathbf{R}}|K(u)|^{m}du\leq m!C^{m} \, .
\end{equation*}%
\medskip

\quad This requirement on the kernel is weaker than the exponential moment
condition imposed by Gao (2003, relation (1.6)).

\quad For each real number $x$, each positive integer $n$ and each $t\in
\lbrack 0,1]$, let us define 
\begin{equation*}
f_{[nt]}(x)=\frac{1}{nh_{n}}\sum_{k=1}^{[nt]}K\Big (\frac{x-X_{k}}{h_{n}}%
\Big )\, ,
\end{equation*}%
where for all $n\geq 1,$ $h_{n}\ $is a strictly positive real number.
Obviously when $t=1$, this is the usual kernel-type estimator of $f$. In
this section we are interested in the moderate deviation principle for the
following processes considered as elements of $D([0,1])$. For fixed real
number $x$, each positive integer $n$ and each $t\in \lbrack 0,1]$, let us
define 
\begin{equation*}
U_{n}(t):=\sqrt{nh_{n}}(f_{[nt]}(x)-\mathbf{E}(f_{[nt]}(x))) \, .
\end{equation*}

\begin{proposition}
\label{thmkerden} Suppose (A.1) and (A.2) hold. Then, assuming that $%
a_{n}\rightarrow 0$ and $a_{n}nh_{n}\rightarrow \infty $, the processes $%
U_{n}(.)$ satisfy (\ref{mdpkey}) with the good rate function $I_{f}(\cdot
)=(f(x)\int K^{2}(u)du)^{-1}I(\cdot )$ where $I(\cdot )$ is defined by (\ref%
{rate}).
\end{proposition}

\textbf{Proof of Proposition \ref{thmkerden}} Let us define 
\begin{equation*}
Y_{n,k}(x)=\frac{1}{\sqrt{h_{n}}}\Big (K\Big (\frac{x-X_{k}}{h_{n}}\Big )-{%
\mathbf{E}}K\Big (\frac{x-X_{k}}{h_{n}}\Big )\Big )\,.
\end{equation*}%
Then 
\begin{equation*}
\frac{\sum_{k=1}^{[nt]}Y_{k,n}(x)}{\sqrt{n}}=\sqrt{nh_{n}}(f_{[nt]}(x)-%
\mathbf{E}(f_{[nt]}(x)))\,.
\end{equation*}%
By stationarity, for any $1\leq j\leq n$, $\sum_{k=1}^{j}\mathrm{Var}%
(Y_{k,n}(x))=j\mathrm{Var}(Y_{1,n}(x))$. Hence the conclusion follows
provided the triangular array of independent centered random variables $%
\{Y_{n,k}(x)\}$ satisfies the conditions of Proposition \ref{forappl}.
Notice that for each $m\geq 1,$ our conditions imply 
\begin{eqnarray*}
{\mathbf{E}}|Y_{n,k}(x)|^{m} &\leq &2^{m}\Big (\frac{1}{h_{n}}\Big )^{m/2}{%
\mathbf{E}}\Big |K\Big (\frac{x-X_{k}}{h_{n}}\Big )\Big |^{m} \\
&\leq &m!2^{m}h_{n}\Vert f\Vert _{\infty }\Big (\frac{1}{h_{n}}\Big )%
^{m/2}C^{m}\,.
\end{eqnarray*}%
Setting $A_{n,k}=2Ch_{n}^{-1/2}$ and $B_{n}=h_{n}\Vert f\Vert _{\infty }$,
the assumptions of Proposition \ref{forappl} hold since $na_{n}h_{n}%
\rightarrow \infty $. It follows that 
\begin{equation*}
\Big \{\frac{\sum_{k=1}^{[nt]}Y_{k,n}(x)}{\sqrt{n}\sqrt{\mathrm{Var}(Y_{n,1})%
}},t\in \lbrack 0,1]\Big \}\text{ satisfies the MDP}.
\end{equation*}%
The proof ends by noticing that the dominated convergence theorem ensures
that 
\begin{equation*}
\mathrm{Var}(Y_{n,1})\rightarrow f(x)\int_{\mathbf{R}}K^{2}(u)du\,.
\end{equation*}

\subsection{Application to a class of dependent variables}

\qquad In the recent years MDP was obtained for classes of dependent random
variables by using various martingale approximation techniques. For example,
papers by Gao (1996), Djellout (2002), Dedecker, Merlev\`{e}de, Peligrad and
Utev (2007), used this approach to obtain MDP for classes of $\phi $-mixing
sequences with polynomial rates. In this section we treat other classes of
mixing sequences by another method: approximating the sums of variables in
blocks with triangular array of independent random variables and applying
then our Theorem \ref{keyresult} to prove the MDP. The measure of
dependence, called $\tau $, that we shall use in this section has been
introduced by Dedecker and Prieur (2004) and it can easily be computed in
many situations such as causal Bernoulli shifts, functions of strong mixing
sequences, iterated random functions and so on. We refer to papers by
Dedecker and Prieur (2004) or Dedecker and Merlev\`{e}de (2006) for precise
estimation of $\tau $ for these examples. Since the rate of convergence in
the next corollary is geometric, we would like also to mention that the
result in this section can be also applied to ARCH models whose coefficients 
$a_{j}$ are zero for large enough $j\geq J$, since these models are
geometrically $\tau $-dependent (see Proposition 5.1 in Comte, Dedecker and
Taupin (2007)).

\medskip

\quad Let us now introduce the dependence coefficients used in what follows.

\medskip

\quad For any real random variable $X$ in ${\mathbf{L}}^{1}$ and any $\sigma 
$-algebra $\mathcal{M}$ of $\mathcal{A}$, let ${\mathbf{P}}_{X|\mathcal{M}}$
be a conditional distribution of $X$ given ${\mathcal{M}}$ and let ${\mathbf{%
P}}_{X}$ be the distribution of $X$. We consider the coefficient $\tau (%
\mathcal{M},X)$ of weak dependence (Dedecker and Prieur, 2004) which is
defined by 
\begin{equation}
\tau (\mathcal{M},X)=\Big \|\sup_{f\in \Lambda _{1}(\mathbf{R})}\Bigr|\int
f(x)\mathbf{P}_{X|\mathcal{M}}(dx)-\int f(x)\mathbf{P}_{X}(dx)\Big |\Big \|%
_{1}\,,  \label{deftau1}
\end{equation}%
where $\Lambda _{1}(\mathbf{R})$ is the set of $1$-Lipschitz functions from $%
\mathbf{R}$ to $\mathbf{R}$.

The $\tau $-coefficient has the following coupling property: If $\Omega $ is
rich enough then the coefficient $\tau (\mathcal{M},X)$ is the infimum of $%
\Vert X-Y\Vert _{1}$ where $Y$ is independent of $\mathcal{M}$ and
distributed as $X$ (see Lemma 5 in Dedecker and Prieur (2004)). This
coupling property allows to relate the $\tau $-coefficient with the strong
mixing coefficient Rosenblatt (1956) defined by 
\begin{equation*}
\alpha (\mathcal{M},\sigma (X))=\sup_{A\in \mathcal{M},B\in \sigma (X)}|{%
\mathbf{P}}(A\cap B)-{\mathbf{P}}(A)P(B)|\,\text{\ .}
\end{equation*}%
as shown in Rio (2000), page 161 (see Peligrad ( 2002) for the unbounded
case). In case when $X$ is bounded, we have 
\begin{equation*}
\tau (\mathcal{M},X)\leq 4\Vert X\Vert _{\infty }\alpha (\mathcal{M},\sigma
(X))\,.
\end{equation*}
For equivalent definitions of the strong mixing coefficient we refer for
instance to Bradley (2007, Lemma 4.3 and Theorem 4.4).

\medskip

\quad If $Y$ is a random variable with values in $\mathbf{R}^{k}$, the
coupling coefficient $\tau $ is defined as follows: If $Y\in {\mathbf{L}}%
^{1}(\mathbf{R}^{k})$, 
\begin{equation}
\tau (\mathcal{M},Y)=\sup \{\tau (\mathcal{M},f(Y)),f\in \Lambda _{1}(%
\mathbf{R}^{k})\}\,,  \label{deftau1}
\end{equation}%
where $\Lambda _{1}(\mathbf{R}^{k})$ is the set of $1$-Lipschitz functions
from $\mathbf{R}^{k}$ to $\mathbf{R}$.

\medskip

\quad We can now define the coefficient $\tau $ for a sequence $%
(X_{i})_{i\in {\mathbf{Z}}}$ of real valued random variables.

\noindent For a strictly sequence $(X_{i})_{i\in {\mathbf{Z}}}$ of
real-valued random variables and for any positive integer $i$, define 
\begin{equation}
\tau (i)=\sup_{k\geq 0}\max_{1\leq \ell \leq k}\frac{1}{\ell }\sup \Big \{%
\tau (\mathcal{M}_{0},(X_{j_{1}},\cdots ,X_{j_{\ell }})),\,i\leq
j_{1}<\cdots <j_{\ell }\Big \}\,,  \label{deftau2}
\end{equation}%
where $\mathcal{M}_{0}=\sigma (X_{j},j\leq 0)$ and supremum also extends for
all $i\leq j_{1}<\cdots <j_{\ell }.$\newline
On an other hand, the sequence of strong mixing coefficients $(\alpha
(i))_{i>0}$ is defined by: 
\begin{equation*}
\ {\alpha }(i)=\alpha (\mathcal{M}_{0},\sigma (X_{j},j\geq i))\,.
\end{equation*}%
In the case where the variables are bounded the following bound is valid 
\begin{equation}
\tau (i)\leq 4\Vert X_{0}\Vert _{\infty }\alpha (i)\,  \label{comptaualpha}
\end{equation}%
(see Lemma 7 Dedecker and Prieur, 2004). \medskip

\quad In the next proposition, we consider a strictly stationary sequence
whose $\tau $-dependence coefficients are geometrically decreasing.

\begin{proposition}
\label{thmtau} Let $(X_{i})_{i\in {\mathbf{Z}}}$ be a strictly stationary
sequence of centered random variables such that $\Vert X_{0}\Vert _{\infty
}<\infty $. Let $S_{n}=\sum_{i=1}^{n}X_{i}$ and $\sigma _{n}^{2}=\mathrm{Var}%
(S_{n})$. Let $(\tau (n))_{n\geq 1}$ be the sequence of dependence
coefficients of $(X_{i})_{i\in {\mathbf{Z}}}$ defined by (\ref{deftau2}).
Assume that $\sigma _{n}^{2}\rightarrow \infty $ and that there exists $\rho
\in ]0,1[$ such that $\tau (n)\leq \rho ^{n}$. Then, for all positive
sequences $a_{n}$ with $a_{n}\rightarrow 0$ and $na_{n}^{2}\rightarrow
\infty $, the normalized partial sums processes $\{\sigma
_{n}^{-1}\sum_{i=1}^{[nt]}X_{i},t\in \lbrack 0,1]\}$ satisfy (\ref{mdpkey})
with the good rate function given in Theorem \ref{keyresult}.
\end{proposition}

\begin{remark}
Taking into account the bound (\ref{comptaualpha}), Corollary \ref{thmtau}
directly applies to strongly mixing sequences of bounded random variables
with geometric mixing rate $({\alpha }(n)\leq \rho ^{n})$.
\end{remark}

\begin{remark}
\label{remvar}Notice that, since the variables are centered and bounded, we
get that $|\mathrm{Cov}(X_{0},X_{k})|\leq \Vert X_{0}\Vert _{\infty }\Vert {%
\mathbf{E}}(X_{k}|\mathcal{M}_{0})\Vert _{1}$. Now from the definition of
the $\tau $-dependence coefficient we clearly have that $\Vert {\mathbf{E}}%
(X_{k}|\mathcal{M}_{0})\Vert _{1}\leq \tau (k)$. It follows that the
condition on the sequence of coefficients $\tau (n)$ implies $\sum_{k}k|%
\mathrm{Cov}(X_{0},X_{k})|<\infty $. This condition together with the fact
that $\sigma _{n}^{2}\rightarrow \infty $ entail that $n^{-1}\mathrm{Var}%
(S_{n})$ converges to a finite number $\sigma ^{2}>0$ (see Lemma 1 in
Bradley (2007)).
\end{remark}

\noindent \textbf{Proof of Proposition \ref{thmtau}.} Let $\varepsilon
_{n}^{2}\rightarrow 0$ in such a way that $\varepsilon
_{n}^{2}na_{n}^{2}\rightarrow \infty $ (this is possible because $%
na_{n}^{2}\rightarrow \infty $) and $\varepsilon _{n}^{2}na_{n}/\log
(a_{n}n)\rightarrow \infty $. Take $p_{n}=\varepsilon _{n}na_{n}$ and $%
q_{n}:=\varepsilon _{n}^{2}na_{n}$.

\quad We now divide the variables $\{X_{i}\}$ in big blocks of size $p_{n}$
and small blocks of size $q_{n}$ in the following way : Let us set $%
\displaystyle k_{n}=[n(p_{n}+q_{n})^{-1}]$. For a given positive integer $n$%
, the set ${1,2,\cdots ,n}$ is being partitioned into blocks of consecutive
integers, the blocks being $I_{1},J_{1},...,I_{k_{n}},J_{k_{n}}$, such that
for each $1\leq j\leq k_{n}$, $I_{j}$ contains $p_{n}$ integers and $J_{j}$
contains $q_{n}$ integers.

\quad Denote by $Y_{j,n}:=\sum_{i\in I_{j}}X_{i}$ and $Z_{j,n}:=\sum_{i\in
J_{j}}X_{i}$ for $1\leq j\leq k_{n}$. Now we consider the following
decomposition: for any $t\in \lbrack 0,1]$, 
\begin{equation}
\sum_{i=1}^{[nt]}X_{i}=\sum_{j=1}^{[k_{n}t]}Y_{j,n}+%
\sum_{j=1}^{[k_{n}t]}Z_{j,n}+R_{n,t}\,,  \label{pip2}
\end{equation}%
where 
\begin{equation*}
R_{n,t}:=\sum_{i=1}^{[nt]}X_{i}-\Big (\sum_{j=1}^{[k_{n}t]}Y_{j,n}+%
\sum_{j=1}^{[k_{n}t]}Z_{j,n}\Big )\,.
\end{equation*}%
\qquad The idea of the proof is the following: Using Lemma 5 in Dedecker and
Prieur (2004), we get the existence of independent random variables $%
(Y_{i,n}^{\ast })_{1\leq i\leq k_{n}}$ with the same distribution as the
random variables $Y_{i,n}$ such that 
\begin{equation}
\mathbf{E}|Y_{i,n}-Y_{i,n}^{\ast }|\leq p_{n}\tau (q_{n})\,.  \label{coupY}
\end{equation}%
Then we show that the partial sums processes $\{\sigma
_{n}^{-1}\sum_{j=1}^{[k_{n}t]}Y_{j,n}^{\ast },t\in \lbrack 0,1]\}$ satisfy (%
\ref{mdpkey}) with the good rate function given in Theorem \ref{keyresult},
while the remainder is negligible for the convergence in distribution, i.e
for all $\eta >0$, 
\begin{equation}
\limsup_{n\rightarrow \infty }a_{n}\log \Big ({\mathbf{P}}\big (\sup_{0\leq
t\leq 1}\frac{\sqrt{a_{n}}}{\sigma _{n}}\big |\sum_{i=1}^{[nt]}X_{i}-%
\sum_{j=1}^{[k_{n}t]}Y_{j,n}^{\ast }\big |\geq \eta \big )\Big )=-\infty \,.
\label{firstapprox}
\end{equation}%
\qquad By stationarity, $\mathrm{Var}(Y_{j,n}^{\ast })=\sigma _{p_{n}}^{2}$
for any $1\leq j\leq k_{n}$ and that, by Remark \ref{remvar}, $k_{n}\sigma
_{p_{n}}^{2}/\sigma _{n}^{2}\rightarrow 1$. Also for any $k\in \lbrack
1,k_{n}]$, $\sum_{j=1}^{k}\mathrm{Var}(Y_{j,n}^{\ast })/(k_{n}\sigma
_{p_{n}}^{2})=k/k_{n}$. Hence, by taking into account these considerations,
we shall verify the conditions of Theorem \ref{keyresult} for the variables $%
\{Y_{j,n}^{\ast }\}_{1\leq j\leq k_{n}}$. According to Comment \ref{comment}
and using stationarity, it suffices to verify that there is a constant $%
C_{1} $ with the property that any $\epsilon >0$ and any $\beta >0$ there is 
$N(\epsilon ,\beta )$ such that for $n>N(\epsilon ,\beta )$ 
\begin{equation}
a_{n}k_{n}{\mathbf{P}}(|S_{p_{n}}|>u\sqrt{a_{n}}\sigma _{n})\leq C_{1}\exp
(-\beta u)\text{\thinspace for any }u\geq \epsilon \,.  \label{condY}
\end{equation}%
Applying Lemma \ref{lmadeddouk} in the Appendix, we derive that there exist
two positive constants $C_{1}$ and $C_{2}$ depending only on $\Vert
X_{0}\Vert _{\infty }$ and $\rho $ such that 
\begin{equation*}
a_{n}k_{n}{\mathbf{P}}(|S_{p_{n}}|>u\sqrt{a_{n}}\sigma _{n})\leq C_{1}\frac{%
a_{n}n}{p_{n}}\exp (-C_{2}\frac{u\sqrt{a_{n}}\sigma _{n}}{\sqrt{p_{n}}})\,.
\end{equation*}%
Since $\sigma _{n}^{2}/n\rightarrow \sigma ^{2}>0$ and by the selection of $%
p_{n}$ we have that $p_{n}=o(na_{n})$ which proves (\ref{condY}) and we
conclude that the process $\{\sigma
_{n}^{-1}\sum_{j=1}^{[k_{n}t]}Y_{j,n}^{\ast },t\in \lbrack 0,1]\}$ satisfies
the conclusion of Theorem \ref{keyresult}.

\quad Hence it remains to show (\ref{firstapprox}). We shall decompose the
proof of this negligibility in several steps.

\quad Using again Lemma 5 in Dedecker and Prieur (2004), there are
independent random variables $(Z_{i,n}^{\ast })_{1\leq i\leq k_{n}}$ with
the same distribution as the random variables $Z_{i,n}$ such that $\mathbf{E}%
|Z_{i,n}-Z_{i,n}^{\ast }|\leq q_{n}\tau (p_{n})$. By the same arguments as
for the sequence $\{Y_{j,n}^{\ast }\}_{1\leq j\leq k_{n}}$, we get that the
Donsker process $\{(k_{n}\sigma
_{q_{n}}^{2})^{-1/2}\sum_{j=1}^{[k_{n}t]}Z_{j,n}^{\ast },t\in \lbrack 0,1]\}$
satisfies (\ref{mdpkey}) with the good rate function given in Theorem \ref%
{keyresult}. Now since $k_{n}\sigma _{q_{n}}^{2}/\sigma _{n}^{2}\sim
q_{n}/p_{n}$ converges to zero as $n\rightarrow \infty $ we easily deduce
that for all $\eta >0$, 
\begin{equation*}
\limsup_{n\rightarrow \infty }a_{n}\log \Big ({\mathbf{P}}\big (\sup_{0\leq
t\leq 1}\frac{\sqrt{a_{n}}}{\sigma _{n}}\big |\sum_{j=1}^{[k_{n}t]}Z_{j,n}^{%
\ast }\big |\geq \eta \big )\Big )=-\infty \,.  \label{secondapprox}
\end{equation*}%
Consequently, to prove (\ref{firstapprox}), it remains to prove that for all 
$\eta >0$, 
\begin{equation}
\limsup_{n\rightarrow \infty }a_{n}\log \Big ({\mathbf{P}}\big (\sup_{0\leq
t\leq 1}\frac{\sqrt{a_{n}}}{\sigma _{n}}\big |%
\sum_{j=1}^{[k_{n}t]}(Y_{j,n}-Y_{j,n}^{\ast }+Z_{j,n}-Z_{j,n}^{\ast })\big |%
\geq \eta \big )\Big )=-\infty \,,  \label{3approx}
\end{equation}%
and 
\begin{equation}
\limsup_{n\rightarrow \infty }a_{n}\log \Big ({\mathbf{P}}\big (\sup_{0\leq
t\leq 1}\frac{\sqrt{a_{n}}|R_{n,t}|}{\sigma _{n}}\geq \eta \big )\Big )%
=-\infty \,.  \label{pip3}
\end{equation}%
By using Markov inequality, we clearly have that 
\begin{eqnarray*}
{\mathbf{P}}\big (\sup_{0\leq t\leq 1}\frac{\sqrt{a_{n}}}{\sigma _{n}}\big |%
\sum_{j=1}^{[k_{n}t]}(Y_{j,n}-Y_{j,n}^{\ast }+Z_{j,n}-Z_{j,n}^{\ast })\big |%
\geq \eta \big ) &\leq &\frac{\sqrt{a_{n}}}{\eta \sigma _{n}}k_{n}({\mathbf{E%
}}|Y_{1,n}-Y_{1,n}^{\ast }|+{\mathbf{E}}|Z_{1,n}-Z_{1,n}^{\ast }|) \\
&\leq &\frac{2}{\eta }\frac{n\sqrt{a_{n}}}{\sigma _{n}}\tau (q_{n})\leq 
\frac{2}{\eta }\frac{n\sqrt{a_{n}}}{\sigma _{n}}e^{-q_{n}\log (1/\rho )}\,,
\end{eqnarray*}%
which proves (\ref{3approx}) by using the selection of $\varepsilon _{n}$
and $q_{n}$ and the fact that $\sigma _{n}^{2}/n\rightarrow \sigma ^{2}>0$.

\quad Since for any $t\in \lbrack 0,1]$, $R_{n,t}$ contains at most $%
2(p_{n}+q_{n})$ terms, by stationarity we have 
\begin{equation*}
{\mathbf{P}}\big (\sup_{0\leq t\leq 1}\frac{\sqrt{a_{n}}|R_{n,t}|}{\sigma
_{n}}\geq \eta \big )\leq (k_{n}+1){\mathbf{P}}\big (\max_{1\leq j\leq
2(p_{n}+q_{n})}\frac{\sqrt{a_{n}}|\sum_{i=1}^{j}X_{i}|}{\sigma _{n}}\geq
\eta \big )\,.
\end{equation*}%
Applying Lemma \ref{lmadeddouk} in the Appendix, we derive that there exist
positive constants $C_{1}$ and $C_{2}$ depending only on $\Vert X_{0}\Vert
_{\infty }$ and $\rho $ such that 
\begin{equation*}
{\mathbf{P}}\big (\max_{1\leq j\leq 2(p_{n}+q_{n})}\frac{\sqrt{a_{n}}%
|\sum_{i=1}^{j}X_{i}|}{\sigma _{n}}\geq \eta \big )\leq C_{1}\exp \Big (%
-C_{2}\frac{\eta \sigma _{n}}{\sqrt{a_{n}}\sqrt{2(p_{n}+q_{n})}}\Big )\,.
\end{equation*}%
It follows that 
\begin{equation*}
a_{n}\log \Big ({\mathbf{P}}\big (\sup_{0\leq t\leq 1}\frac{\sqrt{a_{n}}%
|R_{n,t}|}{\sigma _{n}}\geq \eta \big )\Big )\leq a_{n}\log
(k_{n}+1)+a_{n}\log (C_{1})-C_{2}\frac{\eta \sqrt{a_{n}}\sigma _{n}}{\sqrt{%
2(p_{n}+q_{n})}}\,.
\end{equation*}%
Since $\sigma _{n}^{2}/n\rightarrow \sigma ^{2}>0$ and $p_{n}=o(na_{n})$, we
get that $\sqrt{a_{n}}\sigma _{n}/\sqrt{p_{n}+q_{n}}\rightarrow \infty .$ In
addition $k_{n}\sim a_{n}^{-1}$ implying that $a_{n}\log
(k_{n}+1)\rightarrow 0$. Hence (\ref{pip3}) is proved which completes the
proof of (\ref{firstapprox}) and then of the proposition.

\section{Proofs}

\subsection{\textbf{Proof of Theorem \protect\ref{keyresult}}}

\qquad To prove this theorem we shall use a truncation argument. Without
restricting the generality we shall assume in this proof $%
s_{n}^{2}=\sum_{j=1}^{k_{n}}{\mathbf{E}}(X_{nj})^{2}=1.$ This is possible by
dividing all variables by $s_{n}^{2}$ and redenoting them also by $X_{nj}.$
We truncate the variables in the following way: For all $1\leq j\leq k_{n}$,
let 
\begin{equation*}
X_{nj}^{^{\prime }}:=X_{nj}I(|X_{nj}|\leq \sqrt{a_{n}})-{\mathbf{E}}%
(X_{nj}I(|X_{nj}|\leq \sqrt{a_{n}})\,,
\end{equation*}%
\begin{eqnarray*}
X_{nj}^{"}:= &&X_{nj}I(\sqrt{a_{n}}<|X_{nj}|\leq 1/\sqrt{a_{n}})-{\mathbf{E}}%
(X_{nj}I(\sqrt{a_{n}}<|X_{nj}|\leq 1/\sqrt{a_{n}})) \\
:= &&\bar{X}_{nj}I(|\bar{X}_{nj}|>\sqrt{a_{n}})-{\mathbf{E}}(\bar{X}_{nj}I(|%
\bar{X}_{nj}|>\sqrt{a_{n}})))\,,
\end{eqnarray*}%
and%
\begin{equation*}
X_{nj}^{^{\prime \prime \prime }}:=X_{nj}I(|X_{nj}|>1/\sqrt{a_{n}})-{\mathbf{%
E}}(X_{nj}I(|X_{nj}|>1/\sqrt{a_{n}}))\,.
\end{equation*}%
Above we used also the notation: $\bar{X}_{nj}=X_{nj}I(|X_{nj}|\leq 1/\sqrt{%
a_{n}}).$ Notice first that, since $s_{n}^{2}=1$ 
\begin{equation*}
\sqrt{a_{n}}\sum_{j=1}^{k_{n}}{\mathbf{E}}(|X_{nj}|I(|X_{nj}|>1/\sqrt{a_{n}}%
))\leq a_{n}\rightarrow 0\text{ ,}
\end{equation*}%
and that for any $\delta >0$, 
\begin{equation*}
a_{n}\log {\mathbf{P}}\big (\sum_{j=1}^{k_{n}}|X_{nj}|I(|X_{nj}|>1/\sqrt{%
a_{n}})\geq \delta \big )\leq a_{n}\log {\mathbf{P}}\big (\max_{1\leq j\leq
k_{n}}|X_{nj}|>1/\sqrt{a_{n}}\big )\,.
\end{equation*}%
Hence, by taking into account condition (\ref{neg2}), the variables $%
X_{nj}^{^{\prime \prime \prime }}$ have a negligible contribution to the MDP
(see Theorem 4.2.13 in Dembo and Zeitouni (1998)). Consequently, without
restricting the generality we have just to consider the sums $%
S_{nl}^{^{\prime }}=\sum_{j=1}^{l}X_{nj}^{^{\prime }}$ and $%
S_{nl}^{"}=\sum_{j=1}^{l}X_{nj}^{"}$. We denote by $W_{n}^{^{\prime }}(t)$
(respectively by $W_{n}^{"}(t)$) the random function on $[0,1]$ that is
linear on each interval $[s_{n,i-1}^{2}$, $s_{ni}^{2}]$ and has the values $%
W_{n}^{^{\prime }}(s_{ni}^{2})=S_{ni}^{^{\prime }}$ (respectively $%
W_{n}^{"}(s_{ni}^{2})=S_{ni}^{"}$) at the points of division. Then 
\begin{equation*}
W_{n}(t)\approx W_{n}^{^{\prime }}(t)+W_{n}^{"}(t)\,.
\end{equation*}%
\qquad We first show that the sequence $W_{n}^{"}(t)$ is also negligible,
that is for any $\delta >0$, 
\begin{equation}
\lim_{n\rightarrow \infty }a_{n}\log \left( {\mathbf{P}}\big (\sup_{0\leq
t\leq 1}\sqrt{a_{n}}|W_{n}^{"}(t)|\geq \delta \big )\right) =-\infty \,.
\label{neg1}
\end{equation}

\quad Notice that, $\sqrt{a_{n}}\sum_{k=1}^{k_{n}}|{\mathbf{E}}(\bar{X}%
_{nk}I(|\bar{X}_{nk}|>\sqrt{a_{n}})|$ converges to zero as a consequence of
condition (\ref{lindeberg2}). Hence we have to establish for any $\delta >0$%
\begin{equation*}
\lim_{n\rightarrow \infty }a_{n}\log \Big ({\mathbf{P}}\Big (\sqrt{a_{n}}%
\sum_{k=1}^{k_{n}}\big |\bar{X}_{nk}I(|\bar{X}_{nk}|>\sqrt{a_{n}})\big |\geq
\delta \Big )\Big )=-\infty \,.
\end{equation*}%
Clearly, for any $\lambda >0$ 
\begin{align*}
a_{n}\log \Big ({\mathbf{P}}\Big (\sqrt{a_{n}}\sum_{k=1}^{k_{n}}\big |\bar{X}%
_{nk}I(|\bar{X}_{nk}|& >\sqrt{a_{n}})\big |\geq \delta \Big )\Big )\leq
-\lambda \delta + \\
a_{n}\sum_{k=1}^{k_{n}}\log {\mathbf{E}}\big (\exp (\frac{\lambda |\bar{X}%
_{nk}|}{\sqrt{a_{n}}}I(|\bar{X}_{nk}|& >\sqrt{a_{n}})\big )\,,
\end{align*}%
which shows that it is enough to prove that there is a positive constant $C$
such that for each $\lambda >0$ 
\begin{equation*}
\lim \sup_{n\rightarrow \infty }a_{n}\sum_{k=1}^{\ell _{n}(j)}\log {\mathbf{E%
}}\big (\exp (\frac{\lambda |\bar{X}_{nk}|}{\sqrt{a_{n}}}I(|\bar{X}_{nk}|>%
\sqrt{a_{n}})\big )\leq C\,.
\end{equation*}%
Since $e^{xI(A)}-1=(e^{x}-1)I(A)$ and also $\log (1+x)\leq x$ the above
inequality is implied by

\begin{equation*}
\lim \sup_{n\rightarrow \infty }a_{n}\sum_{k=1}^{k_{n}}{\mathbf{E}}\big (%
\lbrack \exp (\frac{\lambda |\bar{X}_{nk}|}{\sqrt{a_{n}}})-1]I(|\bar{X}%
_{nk}|>\sqrt{a_{n}})\big )\leq C\,,
\end{equation*}%
which is a consequence of condition (\ref{lindebergmain2}).

\medskip

\quad In order to prove that the sequence $W_{n}^{^{\prime }}(t)$ satisfies
the moderate deviation principle, according to Theorem 3.2. in Arcones
(2003-b), it is enough to show that, for a fixed integer $m$, and each $%
0=t_{0}\leq t_{1}\leq \cdots \leq t_{m}\leq 1$, 
\begin{eqnarray}
&&\big (W_{n}^{\prime }(t_{1}),\cdots ,W_{n}^{\prime }(t_{m})\big)\text{
satisfies the MDP in ${\mathbf{R}}^{m}$ with speed $a_{n}$ and the good}
\label{finiteconv} \\
&&\text{ rate function } I_{m}(u_{1},\cdots ,u_{m})=\sum_{\ell =1}^{m}\frac{1%
}{2}\frac{(u_{\ell }-u_{\ell -1})^{2}}{(t_{\ell }-t_{\ell -1})}\text{ with $%
u_{0}=0,$}\,  \notag
\end{eqnarray}%
and for each $\delta >0$ 
\begin{equation}
\lim_{\eta \rightarrow 0}\limsup_{n\rightarrow \infty }a_{n}\log \Big ({%
\mathbf{P}}\Big \{\sup_{|s-t|\leq \eta ,0\leq s,t\leq 1}\sqrt{a_{n}}%
|W_{n}^{\prime }(t)-W_{n}^{\prime }(s)|\geq \delta \Big \}\Big )=-\infty \, .
\label{tight}
\end{equation}%
By the contraction principle (see Theorem 4.2.1 in Dembo and Zeitouni
(1998)) to prove the convergence of the finite dimensional distributions we
have to show that

$Y_{n}^{^{\prime }}:=\big (W_{n}^{^{\prime }}(t_{1}),W_{n}^{^{\prime
}}(t_{2})-W_{n}^{^{\prime }}(t_{1}),\cdots ,W_{n}^{^{\prime
}}(t_{m})-W_{n}^{^{\prime }}(t_{m-1})\big )$ satisfies the MDP in ${\mathbf{R%
}}^{m}$ with speed $a_{n}$ and the good rate function given by 
\begin{equation}
I_{m}^{\prime }(u_{1},\cdots ,u_{m})=\sum_{\ell =1}^{m}\frac{1}{2}\frac{%
u_{\ell }^{2}}{(t_{\ell }-t_{\ell -1})}\,.  \label{iprime}
\end{equation}%
According to Theorem II.2 in Ellis (1984) and independence we have to verify
for each $j,1\leq j\leq m$, 
\begin{equation}
\lim_{n\rightarrow \infty }a_{n}\log \left( {\mathbf{E}}\Big \{\exp \big (%
\frac{1}{\sqrt{a_{n}}}\lambda _{j}(W_{n}^{^{\prime }}(t_{j})-W_{n}^{^{\prime
}}(t_{j-1}))\Big \}\right) =\frac{1}{2}\lambda _{j}^{2}(t_{j}-t_{j-1})\,.
\label{ellisprime}
\end{equation}%
Notice that 
\begin{equation}
W_{n}^{^{\prime }}(t_{j})-W_{n}^{^{\prime }}(t_{j-1})=\sum_{k=\ell
_{n}(t_{j-1})+1}^{\ell _{n}(t_{j})}X_{nk}^{^{\prime }}\,,  \label{decWnprime}
\end{equation}%
with $\ell _{n}(t_{j})$ the maximum $k$ for which $s_{nk}^{2}\leq t_{j}$
(this difference is understood to be 0 if $\ell _{n}(t_{j-1})=\ell
_{n}(t_{j}))$). We shall verify the conditions of Lemma 2.3 in Arcones
(2003-a), given for convenience in Appendix, to the real valued random sum
of independent random variables: $Y_{n1}=X_{n,\ell _{n}(t_{j-1})+1}^{\prime
},\dots ,Y_{nk_{n}}=X_{n,\ell _{n}(t_{j})}^{\prime }$. Since the random
variables $X_{nk}^{\prime }$ are uniformly bounded by $\sqrt{a_{n}}$,
condition (\ref{bound}) holds. Now the Lindeberg's condition (\ref%
{lindeberg2}) clearly implies (\ref{lind}). Hence it remains to verify that
for any $1\leq j\leq m$, 
\begin{equation}
\lim_{n\rightarrow \infty }\sum_{k=\ell _{n}(t_{j-1})+1}^{\ell _{n}(t_{j})}{%
\mathbf{E}}(X_{nk}^{^{\prime }})^{2}=(t_{j}-t_{j-1})\,.  \label{p1tay}
\end{equation}%
By condition (\ref{lindeberg2}), (\ref{p1tay}) holds provided that 
\begin{equation}
\lim_{n\rightarrow \infty }\sum_{k=\ell _{n}(t_{j-1})+1}^{\ell _{n}(t_{j})}{%
\mathbf{E}}(X_{nk})^{2}=(t_{j}-t_{j-1})\,.  \label{p3tay}
\end{equation}%
For $n$ sufficiently large $\ell _{n}(t_{j-1})\neq \ell _{n}(t_{j})$ and $%
\sum_{k=\ell _{n}(t_{j-1})+1}^{\ell _{n}(t_{j})}{\mathbf{E}}%
(X_{nk})^{2}=s_{n,\ell _{n}(j)}^{2}-s_{n,\ell _{n}(j-1)}^{2}$. Also,
condition (\ref{lindeberg2}) implies that for all $k$, $s_{nk}^{2}%
\rightarrow 0,$ therefore 
\begin{equation}
\lim_{n\rightarrow \infty }s_{n,\ell _{n}(t_{j})}^{2}=t_{j}\,.
\label{conslindsnj}
\end{equation}%
Hence (\ref{p3tay}) holds, and so does (\ref{p1tay}). This ends the proof of
(\ref{ellisprime}) and of the (\ref{finiteconv}).

\medskip

\quad To prove (\ref{tight}), we notice that by Theorem 7.4 in Billingsley
(1999), for each $\delta >0$, 
\begin{equation*}
{\mathbf{P}}\big (\sup_{d(s,t)\leq \eta }\sqrt{a_{n}}|W_{n}^{\prime
}(t)-W_{n}^{\prime }(s)|\geq 3\delta \big )\leq \sum_{i=1}^{m}{\mathbf{P}}%
\big (\sup_{\frac{i-1}{m}\leq s<\frac{i}{m}}\sqrt{a_{n}}|W_{n}^{\prime
}(s)-W_{n}^{\prime }(\frac{i-1}{m})|\geq \delta \big )\,,
\end{equation*}%
where $m=[\delta ^{-1}]$. In terms of partial sums and above notation, we
get 
\begin{equation*}
{\mathbf{P}}\big (\sup_{\frac{i-1}{m}\leq s<\frac{i}{m}}\sqrt{a_{n}}%
|W_{n}^{\prime }(s)-W_{n}^{\prime }(\frac{i-1}{m})|\geq \delta \big )\leq {%
\mathbf{P}}\big (\max_{\ell (\frac{i-1}{m})\leq k\leq \ell (\frac{i}{m})+1}%
\sqrt{a_{n}}|\sum_{j=\ell (\frac{i-1}{m})}^{k}X_{nj}^{\prime })|\geq \delta %
\big )\,.
\end{equation*}%
By Lindeberg's condition (\ref{lindeberg2}) and (\ref{conslindsnj}) and with
the notation $B_{i,m}^{2}=\sum_{j=\ell (\frac{i-1}{m})+1}^{\ell (\frac{i}{m}%
)}{\mathbf{E}}(X_{nj}^{\prime })^{2}$ we have 
\begin{equation}
\lim_{n\rightarrow \infty }B_{i,m}^{2}=\lim_{n\rightarrow \infty
}\sum_{j=\ell (\frac{i-1}{m})}^{\ell (\frac{i}{m})}B_{i,m}^{2}=\frac{1}{m}
\label{rel}
\end{equation}%
This limit along with Kolmogorov maximal inequality (2.13) in Petrov (1995)
and the fact that $a_{n}\rightarrow 0$ as $n\rightarrow \infty $ gives 
\begin{equation*}
{\mathbf{P}}\big (\max_{\ell (\frac{i-1}{m})\leq k\leq \ell (\frac{i}{m})+1}%
\sqrt{a_{n}}|\sum_{j=\ell (\frac{i-1}{m})}^{k}X_{nj}^{\prime }|\geq \delta %
\big )\leq 2{\mathbf{P}}(\sqrt{a_{n}}|\sum_{j=\ell (\frac{i-1}{m})\ }^{\ell (%
\frac{i}{m})+1}X_{nj}^{\prime }|\geq 2^{-1}\delta )
\end{equation*}%
for all $n$ sufficiently large. Now we apply Prokhorov's inequality (see
Lemma \ref{lmaprok} in the Appendix) with $B=2\sqrt{a_{n}}$ and $t=\frac{1}{2%
}\delta a_{n}^{-1/2}$ to obtain 
\begin{equation*}
{\mathbf{P}}(\sqrt{a_{n}}|\sum_{j=\ell (\frac{i-1}{m})}^{\ell (\frac{i}{m}%
)+1}X_{nj}^{^{\prime }}|\geq 2^{-1}\delta )\leq \exp (-\frac{\delta }{4a_{n}}%
\text{ arcsinh}\frac{\delta }{2B_{i,m}^{2}})\,.
\end{equation*}%
Therefore by (\ref{rel}), 
\begin{equation*}
\lim \sup_{n\rightarrow \infty }a_{n}\log {\mathbf{P}}(\sqrt{a_{n}}%
|\sum_{j=\ell (\frac{i-1}{m})}^{\ell (\frac{i}{m})+1}X_{nj}^{^{\prime
}}|\geq 2^{-1}\delta )\leq -\frac{\delta }{4}\text{ arcsinh}\frac{m\delta }{2%
}\,,
\end{equation*}%
which converges to $-\infty $ when $m\rightarrow \infty $. This convergence
implies (\ref{tight}).

\subsection{Proof of Corollary \protect\ref{onenonbis}}

\qquad We just have to prove that Condition (\ref{onecondm}) implies
Condition (\ref{lindebergmain2}) (condition (\ref{neg2}) being obviously
satisfied). The proof is straightforward but delicate and it is inspired by
the type of arguments developed by Arcones (2003-a, Theorem 2.4) and
Djellout (2002). So, according to Comment \ref{comment}, we shall verify
that condition (\ref{vv}) holds for all $1\leq u\leq 1/a_{n}$.

\quad Recall that $f(x)$ and $g(x)$ are strictly increasing continuous
functions and for any positive integer $m$, $\ f(m)=a_{m}s_{m}^{2}$ $\ $and$%
\ \ g(m)=s_{m}^{2}/a_{m}$. Fix an integer $n_{0}$ and for any $n\geq n_{0}$
and $1\leq u\leq 1/a_{n}$ define $N=N(u,n)$ as $N=g^{-1}(u^{2}f(n))$. So $%
n=f^{-1}(u^{-2}g(N)).$ Notice that $N$ might not be an integer. Obviously,
by monotonicity 
\begin{equation*}
g^{-1}(\ f(n))\leq N\leq g^{-1}(\ f(n)/a_{n}^{2})=n
\end{equation*}%
and 
\begin{equation*}
N\leq n\leq f^{-1}(g(N))\,.
\end{equation*}%
Notice that by the above relations and the assumption $a_{n}s_{n}^{2}%
\rightarrow \infty $, $N(n)\rightarrow \infty $ as $n\rightarrow \infty $
uniformly in $u\geq 1$.

\quad Now, by the definitions of $f(x)$, $g(x),$ and $N,$\ \ for any $1\leq
i\leq k_{n}$ we have%
\begin{equation*}
P(|X_{n,i}|\geq u\sqrt{a_{n}}s_{n})\leq \sup_{N\leq n\leq
f^{-1}(g(N))}\sup_{1\leq i\leq k_{n}}P(|X_{n,i}|\geq s_{[N]}/\sqrt{a_{[N]}})
\, .
\end{equation*}%
Whence, by condition (\ref{onecondm}) for any $\lambda >0,$ and $n\geq
n_{0}(\lambda )$

\begin{equation*}
k_{[N]}P(|X_{n,i}|\geq u\sqrt{a_{n}}s_{n})\leq \exp (-\lambda /a_{[N]})\,.
\end{equation*}
So, for $u\geq 1$ and $n\geq n_{0}(\lambda )$ and taking into account that
both $s_{n}^{2}a_{n}$ and $s_{n}^{2}/k_{n}$ are nondecreasing along with
condition (\ref{lindeberg2}), we easily derive the sequence of inequalities 
\begin{gather*}
a_{n}\sum_{i=1}^{k_{n}}P(|X_{n,i}|\geq u\sqrt{a_{n}}s_{n})\leq a_{n}k_{n}%
\frac{1}{k_{[N]}}\exp (-\lambda /a_{[N]}) \\
=\frac{f(n)}{s_{n}^{2}}\frac{k_{n}}{k_{[N]}}\exp (-\lambda /a_{[N]})\leq 
\frac{g(N)}{s_{n}^{2}}\frac{k_{n}}{k_{[N]}}\exp (-\lambda /a_{[N]}) \\
\leq \frac{s_{[N]+1}^{2}}{s_{n}^{2}a_{[N]+1}}\frac{k_{n}}{k_{[N]}\ }\exp
(-\lambda /a_{[N]})\leq 2\exp ((-\lambda +1)/a_{[N]}) \, .
\end{gather*}

\quad Hence, it remains to compare $u$ with $1/a_{[N]}$. Recall that $%
u^{2}=g(N)/f(n)$ and by monotonicity of $g(x)$ and $f(n)$ and condition (\ref%
{lindeberg2}), we have 
\begin{equation*}
u\leq \frac{1}{s_{n}\sqrt{a_{n}}}\frac{s_{[N]+1}}{\sqrt{a_{[N]+1}}}\leq 
\frac{1}{s_{n}\sqrt{a_{n}}}\frac{s_{[N]+1}^{2}}{\sqrt{s_{[N]}^{2}a_{[N]}}}%
\leq 2\frac{1}{a_{[N]}}\frac{s_{[N]}\sqrt{a_{[N]}}}{s_{n}\sqrt{a_{n}}}\text{
.}
\end{equation*}%
Since $[N]\leq n$, we get $u\leq 2/a_{[N]}$ and therefore condition (\ref{vv}%
) holds for all $1\leq u\leq 1/a_{n}$.

\subsection{Proof of Proposition \protect\ref{forappl}}

\qquad We just have to verify that condition (\ref{lindmain21}) holds. By
Maclaurin expansion and condition (\ref{mom}) we get that for $n\geq n_{0}\ $
\begin{eqnarray*}
& & a_{n}\sum_{j=1}^{k_{n}}{\mathbf{E}}([\exp \frac{\beta |X_{nj}|}{\sqrt{%
a_{n}}s_{n}}]I(\frac{|X_{nj}|}{\sqrt{a_{n}}s_{n}}>\epsilon )\leq \frac{1}{%
\epsilon ^{3}\sqrt{a_{n}}s_{n}^{3}}\sum_{j=1}^{k_{n}}\sum_{p=0}^{\infty }%
\frac{1}{p!}E\Big ( \frac{\beta ^{p}|X_{nj}|^{3+p}}{(\sqrt{a_{n}}s_{n})^{p}} %
\Big ) \\
& & \leq \frac{24B_{n}}{\epsilon ^{3}\sqrt{a_{n}}s_{n}^{3}}%
\sum_{j=1}^{k_{n}}\sum_{p=0}^{\infty }p^{3}|A_{nj}|^{3} \Big ( \beta \frac{%
|A_{nj}|}{\sqrt{a_{n}}s_{n}} \Big )^{p} \leq \frac{24}{\beta \epsilon ^{3}}%
\frac{B_{n}}{s_{n}^{2}}\sum_{j=1}^{k_{n}}|A_{nj}|^{2}\sum_{p=0}^{\infty } %
\Big ( \frac{8\beta|A_{nj}|}{\sqrt{a_{n}}s_{n}} \Big )^{p+1}
\end{eqnarray*}%
which converges to $0$ by (\ref{condAk1}) together with (\ref{condAk2}).

\section{Appendix}

\quad We first state Lemma 2.3 in Arcones (2003a).

\begin{lemma}
Assume that $(Y_{n1},Y_{n2},\dots ,Y_{nk_{n}})$ is a triangular array of
independent random variables, with mean zero and such that 
\begin{equation*}
\lim_{n\rightarrow \infty }{\mathbf{E}}(\sum_{j=1}^{k_{n}}Y_{nj}^{2})=\sigma
^{2}
\end{equation*}%
Let $(a_{n})_{n\geq 1}$ be a sequence of real numbers converging to zero.
Assume there is a constant $\tau $ such that 
\begin{equation}
\sup_{1\leq i\leq k_{n}}|Y_{nj}|\leq \tau \sqrt{a_{n}}\text{ a.s.}
\label{bound}
\end{equation}%
and for each $\delta >0$%
\begin{equation}
a_{n}\sum_{j=1}^{k_{n}}{\mathbf{P}}(|Y_{nj}|\geq \delta \sqrt{a_{n}}%
)\rightarrow 0\text{ as }n\rightarrow \infty \,.  \label{lind}
\end{equation}%
Then, for any $t\in {\mathbf{R,}}$ 
\begin{equation*}
\ a_{n}\log {\mathbf{E}}\exp \Big (t\frac{\sum_{j=1}^{k_{n}}Y_{nj}}{\sqrt{%
a_{n}}}\Big )\rightarrow \frac{t^{2}\sigma ^{2}}{2}\text{ as }n\rightarrow
\infty
\end{equation*}%
and therefore $\{\sum_{j=1}^{k_{n}}Y_{nj}\}$ satisfies the MDP in with speed 
$a_{n}$ and rate function $I(t)=\frac{t^{2}}{2\sigma ^{2}}$.
\end{lemma}

\quad Now we give the following consequence of Corollary 3 in Dedecker and
Doukhan (2003).

\begin{lemma}
\label{lmadeddouk} Let $(X_{i})_{i\in {\mathbf{Z}}}$ be a strictly
stationary sequence of centered real random variables such that $\Vert
X_{0}\Vert _{\infty }<\infty $. Let $(\tau (n))_{n\geq 1}$ be the sequence
of dependence coefficients of $(X_{i})_{i\in {\mathbf{Z}}}$ defined by (\ref%
{deftau2}). Assume that there exist $\rho \in ]0,1[$ such that $\tau (n)\leq
\rho ^{n}$. Let $S_{k}=\sum_{i=1}^{k}X_{i}$. Then there exist constants $%
C_{1}$ and $C_{2}$ depending only on $\rho $ and $\Vert X_{0}\Vert _{\infty
} $ such that the following inequality holds for any integer $m\geq 1$: 
\begin{equation*}
{\mathbf{P}}(\max_{1\leq j\leq m}|S_{j}|>x)\leq C_{1}\exp (-C_{2}x/\sqrt{m}%
)\,.
\end{equation*}
\end{lemma}

\noindent \textbf{Proof of Lemma \ref{lmadeddouk}} First we notice that by
the definition of the $\tau -$dependence coefficient 
\begin{equation*}
\gamma (n)=\Vert {\mathbf{E}}(X_{n}|{\mathcal{M}}_{0})\Vert _{1}\leq \tau
(n)\,\leq \rho ^{n}.
\end{equation*}%
By stationarity and applying Corollary 3 in Dedecker and Doukhan (2003), we
get that for any $1\leq i\leq j\leq m$, there exists a constant $K$
depending only on $\rho $ and $\Vert X_{0}\Vert _{\infty }$ such that 
\begin{equation*}
{\mathbf{P}}\big (|\sum_{\ell =i}^{j}X_{\ell }|>x\big )\leq K\exp \Big (%
\frac{-x\sqrt{\log (1/\rho )}}{e\Vert X_{0}\Vert _{\infty }\sqrt{j-i+1}}\Big
)\,.
\end{equation*}%
Hence the lemma follows by taking into account Theorem 2.2 in M\'{o}ricz,
Serfling and Stout (1982) together with the remark (ii) stated page 1033 in
their paper.

\medskip

\quad Now we recall the Prokhorov's inequality (1959) that we used in the
paper.

\begin{lemma}
\label{lmaprok} Assume that we have an independent random vector (not
necessarily Stationary) $(X_{1},$ $X_{2},...X_{m})$, centered such that 
\begin{equation*}
\max_{1\leq i\leq m}|X_{i}|\leq B\text{ a.s.}
\end{equation*}%
Denote by $s_{n}^{2}=\sum_{j=1}^{n}E(X_{j}^{2})$ . Then for all $t>0$ , the
following inequality holds 
\begin{equation*}
{\mathbf{P}}(|\sum_{j=1}^{n}X_{j}|\geq t)\leq \exp \big (-\frac{t}{2B}\ 
\text{arcsinh}\frac{Bt}{2s_{n}^{2}}\big )\,.
\end{equation*}
\end{lemma}

\medskip

\quad We turn now to the proof of the Comment \ref{comment}.

\medskip

\noindent \textbf{Proof of Comment \ref{comment}}.

\quad Denote $\bar{X}_{nj}=X_{nj}I[|X_{nj}|<s_{n}/\sqrt{a_{n}}].$ We show
that (\ref{lindebergmain2}) is equivalent to the following condition: There
is a constant $C_{1}$ with the property: for any $\beta >0$ there is $%
N(\beta )$ such that for $n>N(\beta )$ 
\begin{equation}
a_{n}\sum_{j=1}^{k_{n}}{\mathbf{P}}(|\bar{X}_{nj}|>u\sqrt{a_{n}}s_{n})\leq
C_{1}\exp (-\beta u)\text{ for all }u\geq 1\,,  \label{equi}
\end{equation}%
which is equivalent to (\ref{vv}).\newline
For any $u\geq 1$ 
\begin{equation*}
a_{n}\sum_{j=1}^{k_{n}}{\mathbf{P}}(|\bar{X}_{nj}|>u\sqrt{a_{n}}s_{n})\leq
a_{n}\sum_{j=1}^{k_{n}}\exp (-\beta u){\mathbf{E}}(\exp (\beta \frac{|\bar{X}%
_{nj}|}{\sqrt{a_{n}}s_{n}})I(|\bar{X}_{nj}|>u\sqrt{a_{n}}s_{n})\,.
\end{equation*}%
Hence (\ref{lindebergmain2}) implies (\ref{equi}). On the other hand, 
\begin{equation*}
{\mathbf{E}}([\exp \frac{\beta |\bar{X}_{nj}|}{2\sqrt{a_{n}}s_{n}}]I(|\bar{X}%
_{nj}|>\sqrt{a_{n}}s_{n})=e^{\beta /2}{\mathbf{P}}(|\bar{X}_{nj}|>\sqrt{a_{n}%
}s_{n})+\frac{\beta }{2}\int_{1}^{\infty }e^{\beta u/2}{\mathbf{P}}(|\bar{X}%
_{nj}|>u\sqrt{a_{n}}s_{n})du
\end{equation*}%
Now if (\ref{equi}) holds then for $n$ sufficiently large, 
\begin{equation*}
a_{n}\sum_{j=1}^{k_{n}}{\mathbf{E}}([\exp \beta \frac{|\bar{X}_{nj}|}{2\sqrt{%
a_{n}}s_{n}}]I(|\bar{X}_{nj}|>\sqrt{a_{n}}s_{n})\leq C_{1}(e^{-\beta /2}+%
\frac{\beta }{2}\int_{1}^{\infty }e^{-\beta u/2}du)\leq 2C_{1}e^{-\beta
/2}\,,
\end{equation*}%
proving that (\ref{lindebergmain2}) is satisfied.


\begin{thebibliography}{99}
\bibitem{Ar} Arcones, M.A. (2003-a). Moderate deviations of empirical
processes. \textit{Stochastic inequalities and applications.} Progr. Probab. 
\textbf{56}, Birkh\"{a}user, Basel, 189-212.

\bibitem{Ar} Arcones, M.A. (2003-b). The large deviation principle for
stochastic processes I. \textit{Theory of Probability and its Applications}. 
\textbf{47}, 567-583.

\bibitem{Ar} Arcones, M.A. (2003-c). The large deviation principle for
stochastic processes II. \textit{Theory of Probability and its Applications}%
. \textbf{48}, 19-44.

\bibitem{Bi} Billingsley, P. (1999). \textit{Convergence of Probability
Measures}. Wiley, New York.

\bibitem{} Bradley, R.C. (1997). On quantiles and the central limit question
for strongly mixing sequences. \textit{J. Theor. Probab.} \textbf{10},
507-555.

\bibitem{Brad} Bradley, R.C. (2007). \textit{Introduction to strong mixing
conditions.} Vol. 1,2,3. Kendrick Press.

\bibitem{CDT} Comte, F., Dedecker, J. and Taupin, M.L. (2007). Adaptive
density estimation for general ARCH models. \textit{%
http://www.math-info.univ-paris5.fr/~comte/publi.html}

\bibitem{DeDo} Dedecker, J. and Doukhan, P. (2003). A new covariance
inequality and applications. \textit{Stoch. Processes Appl.} \textbf{106,}
63-80.

\bibitem{DMPU} Dedecker, J., Merlev\`ede, F., Peligrad, M. and Utev, S.
(2007). Moderate deviations for stationary sequences of bounded random
variables. \textit{Pr\'epublication 1183. LPMA. Universit\'e Paris 6.}

\bibitem{DePr04} Dedecker, J. and Prieur, C. (2004). Coupling for $\tau $%
-dependent sequences and applications. \textit{J. Theoret. Probab.} \textbf{%
17}, 861--885.

\bibitem{DMV} Dedecker, J. and Merlev\`{e}de F. (2006). Inequalities for
partial sums of Hilbert-valued dependent sequences and applications. \textit{%
Math. Methods Statist.}, \textbf{15}, 176--206.

\bibitem{DZ} Dembo, A. and Zeitouni, O. (1998). \textit{Large Deviations
Techniques and Applications}, 2nd edition. Springer New York.

\bibitem{D} Djellout, H. (2002). Moderate Deviations for Martingale
Differences and applications to $\phi-$mixing sequences. \textit{Stoch.
Stoch. Rep.} \textbf{73}, No.1-2, 37-63.

\bibitem{El} Ellis, R. S. (1984). Large deviations for a general class of
random vectors. \textit{Ann. Probab.} \textbf{12}, 1-12.

\bibitem{FU} Gao, F-Q. (1996). Moderate deviations for martingales and
mixing random processes. \textit{Stochastic Process. Appl.}, \textbf{61},
263--275.

\bibitem{g} Gao, F-Q. (2003). Moderate deviations and large deviations for
kernel density estimators. \textit{J. Theoret. Probab.}, \textbf{16},
401-418.

\bibitem{Led} Ledoux, M. (1992). Sur les d\'{e}viations mod\'{e}r\'{e}es de
sommes de variable al\'{e}atoires vectorielles ind\'{e}pendantes de m\^{e}me
loi. \textit{Ann. Inst. Henri Poincar\'{e}}. \textbf{28}, 267-280.

\bibitem{LeTa} Ledoux, M. and Talagrand, M. (1991). \textit{Probability in
Banach spaces.} Springer-Verlag, Berlin.

\bibitem{l} Louani, D. (1998). Large deviations limit theorems for kernel
density estimator. \textit{Scand.J. Statist.} \textbf{25}, 243-253.


\bibitem{MSS} M\'{o}ricz, F.A., Serfling, R.J. and Stout, W.F. (1982).
Moment and Probability bounds with quasi-superadditive structure for the
maximum partial sum. \textit{Ann. Probab.}, \textbf{10}, 1032-1040.

\bibitem{mpw} Mokkadem, A., Pelletier, M. and Worms, J. (2005). Large and
moderate deviations principles for kernel estimators of a multivariate
density and its partial derivatives. \textit{Aust.N.Z.J.Stat.} \textbf{47},
489-502.

\bibitem{mpw} Mokkadem, A., Pelletier, M. and Thiam, B. (2007). Large and
moderate deviations principles for kernel estimators of a multivariate
regression. \textit{hal-00136115}. \textit{http://arxiv.org/abs/math/0703341}%
.

\bibitem{PEL} Peligrad, M. (2002). Some remarks on coupling of dependent
random variables. \textit{Stat. and Prob. Lett. }\textbf{60}, 201-209.

\bibitem{p} Petrov, V. (1995). \textit{Limit theorems in probability theory}%
. Oxford Studies in Probability Series. Clarendon press Oxford.

\bibitem{P} Puhalskii, A. (1994). Large deviations of semimartingales via
convergence of the predictable characteristics. \textit{Stoch. Stoch. Rep. }%
\textbf{49}, 27-85.

\bibitem{Prokho} Prokhorov, Yu. V. (1959). An extremal problem in
probability theory. \textit{Theory of Probability and its Applications.} 
\textbf{4}, 201-203.

\bibitem{rio} Rio, E. (2000). \textit{Th\'eorie asymptotique des processus
al\'eatoires faiblement d\'ependants}. Math\'ematiques \& Application,%
\textit{\ }\textbf{31}, Springer-Verlag Berlin Heidelberg .

\bibitem{Ro} Rosenblatt, M. (1956). A central limit theorem and a strong
mixing condition. \textit{Proc. Nat. Acad. Sci. U. S. A. } \textbf{42,}
43-47.
\end{thebibliography}
\end{document}